\documentclass[11pt]{article}%
\usepackage{amsmath}%
\usepackage{amsfonts}%
\usepackage{amssymb}%
\usepackage{graphicx}
\providecommand{\U}[1]{\protect\rule{.1in}{.1in}}
\newtheorem{theorem}{Theorem}[section]

\newtheorem{definition}[theorem]{Definition}

\numberwithin{equation}{section}

\begin{document}

\title{Method for solving inverse spectral problems on quantum star graphs}
\author{Sergei A. Avdonin$^{1}$, Vladislav V. Kravchenko$^{2}$\\{\small $^{1}$ Department of Mathematics and Statistics, University of Alaska,
Fairbanks, AK 99775, USA}\\{\small $^{2}$ Departamento de Matem\'{a}ticas, Cinvestav, Unidad
Quer\'{e}taro, }\\{\small Libramiento Norponiente \#2000, Fracc. Real de Juriquilla,
Quer\'{e}taro, Qro., 76230 MEXICO}\\{\small e-mail: s.avdonin@alaska.edu, vkravchenko@math.cinvestav.edu.mx}}
\maketitle

\begin{abstract}
A new method for solving inverse spectral problems on quantum star graphs is
proposed. The method is based on Neumann series of Bessel functions
representations for solutions of Sturm-Liouville equations. The
representations admit estimates for the series remainders which are
independent of the real part of the square root of the spectral parameter.
This feature makes them especially useful for solving direct and inverse
spectral problems requiring calculation of solutions on large intervals in the
spectral parameter. Moreover, the first coefficient of the representation is
sufficient for the recovery of the potential.

The method for solving the inverse spectral problem on the graph consists in
reducing the problem to a two-spectra inverse Sturm-Liouville problem on each
edge. Then a system of linear algebraic equations is derived for computing the
first coefficient of the series representation for the solution on each edge
and hence for recovering the potential. The proposed method leads to an
efficient numerical algorithm that is illustrated by a number of numerical tests.

\end{abstract}

\section{Introduction}

Under quantum graphs or, in other words, differential equation networks we
understand differential operators on metric graphs coupled by certain vertex
matching conditions. Network-like structures play a fundamental role in many
problems of science and engineering. The classical problem here that comes
from applications is the problem of oscillations of the flexible structures
made of strings, beams, cables, and struts. These models describe bridges,
space-structures, antennas, transmission-line posts and other typical objects
of civil engineering. Quantum graphs arise also as natural models of various
phemonema in chemistry (free-electron theory of conjugated molecules), biology
(genetic networks, dendritic trees), geophysics, environmental science,
decease control, and even in Internet (Internet or network tomography). In
physics, interest to quantum graphs arose, in particular, from applications to
nano-electronics and quantum waveguides. On the other hand, quantum graph
theory gives rise to numerous challenging problems related to many areas of
mathematics from combinatorics to PDE and spectral theory. A number of surveys
and collections of papers on quantum graphs appeared last years, and the first
books on this topic by Berkolaiko and Kuchment \cite{BerkolaikoKuchment},
Mugnolo \cite{Mugnolo} and Kurasov \cite{Kurasov} contain excellent lists of references.

Inverse spectral theory of network-like structures is an important part of the
rapidly developing area of applied mathematics --- analysis on graphs. The
known results in this direction concern almost exclusively trees, i.e. graphs
without cycles, see e.g. \cite{Kurasov et al 2005, Yurko2005, Belishev
Vakulenko 2006, AvdoninKurasov2008, AvdoninLeugeringMikhaylov2010,
AvdoninBell2015}.

To date, there are few papers containing numerical results for inverse
problems on graphs, all of them concern only very simple trees \cite{Belishev
Vakulenko 2006, AvdoninBelinskiyMatthews}. It is known that the problems of
space discretization of differential equations on metric graphs turn out to be
very difficult, and even the forward boundary value problems on graphs contain
a lot of numerical challenges (see, e.g. \cite{ArioliBenzi2018}).

In the present work we develop a new method for solving inverse spectral
problems on compact quantum star graphs that leads to efficient numerical
algorithms. This method is based on Neumann series of Bessel functions
representations for solutions of Sturm-Liouville equations. The method was
proposed in \cite{KNT} and successfully applied to solving direct and inverse
spectral problems on intervals and the semiaxis in \cite{KrBook2020}. Here we
extend this method to the inverse spectral problem for the Schr\"odinger
equation on star graphs.

The Neumann series of Bessel functions representations for solutions of
Sturm-Liouville equations,
\[
-y^{\prime\prime}+q(x)y=\rho^{2}y,
\]
possess two remarkable features which make them especially convenient for
solving inverse problems. The remainders of the series admit $\rho
$-independent bounds for $\rho\in\mathbb{R}$, and the potential of the
equation can be recovered from the very first coefficient of the series. The
first feature allows us to work with the approximate solutions on very large
intervals in $\rho$, and the second implies that computationally satisfactory
results require considering a reduced number of the terms of the series which
eventually results in a reduced number of linear algebraic equations which
should be solved in each step.

The method consists of the following steps. First, using the given spectral
data we compute coefficients of the series representations for solutions
satisfying the homogeneous Dirichlet condition at the boundary vertices and
for their derivatives at the end point of each edge, which is associated with
the common vertex of the star graph. This first step allows us to split the
problem on the graph into separate problems on each edge. Second, the set of
coefficients for the series representation of the solution, evaluated at the
end point of the interval allows us to compute the Dirichlet-Dirichlet
eigenvalues of the potential on each edge. Moreover, the first feature of the
series representations implies that if necessary hundreds of the eigenvalues
can be computed with uniform accuracy, and for this few coefficients of the
series representations are sufficient. Next, the knowledge of the
Dirichlet-Dirichlet eigenvalues allows us to compute the important parameter
(half of the integral of the potential) that is necessary for computing the
Dirichlet-Neumann eigenvalues, for which, besides the parameter, the
coefficients from the series representation of the derivative of the solution
are used.

Thus, the inverse spectral problem on a graph is reduced to a two spectra
inverse Sturm-Liouville problem on each edge. Results on the uniqueness and
solvability of the two spectra problem are well known and can be found, e.g.,
in \cite{Chadan et al 1997}, \cite{LevitanInverse}, \cite{SavchukShkalikov},
\cite{Yurko2007}. For this problem we propose a method which again involves
the Neumann series of Bessel functions representations for solutions and their
derivatives. It allows us to compute multiplier constants \cite{Brown et al
2003} relating the Dirichlet-Neumann eigenfunctions associated with the same
eigenvalues but normalized at the opposite endpoints of the interval. This
leads to a system of linear algebraic equations for the coefficients of their
series representations already for interior points of the interval. Solving
the system we find the very first coefficient, from which the potential is
recovered. Note that the Neumann series of Bessel functions representations
were first used for solving inverse Sturm-Liouville problems in
\cite{Kr2019JIIP}. Later on the approach from \cite{Kr2019JIIP} was improved
in \cite{KrBook2020} and \cite{KT2021 IP1}. In these papers the system of
linear algebraic equations was obtained from the Gelfand-Levitan integral
equation. Here we develop another approach which is based on the consideration
of the eigenfunctions normalized at the opposite endpoints and does not
involve the Gelfand-Levitan equation.

\section{Problem setting}

Let $\Omega$ denote a compact star graph consisting of $M$ edges $e_{1}%
$,...,$e_{M}$ connected at the vertex $v$. Every edge $e_{j}$, is identified
with an interval $(0,L_{j})$ of the real line in such a way that zero
corresponds to the boundary vertex $\gamma_{j}$, and the endpoint $L_{j}$
corresponds to the vertex $v$. By $\Gamma$ we denote the set of the boundary
vertices of $\Omega$, $\Gamma=\left\{  \gamma_{1},\ldots,\gamma_{M}\right\}
$. We associate the standard spectral problem on the graph $\Omega$ with a
real valued potential $q\in L_{2}(\Omega)$:
\begin{equation}
-w^{\prime\prime}(x)+q(x)w(x)=\lambda w(x), \label{Schr}%
\end{equation}%
\begin{equation}
w\in C(\Omega), \label{wcontinuity}%
\end{equation}%
\begin{equation}
\sum_{j=1}^{M} \partial w_{j}(v)=0, \label{KN}%
\end{equation}%
\begin{equation}
w(\gamma_{j})=0,\quad j=1,\ldots,M, \label{Dirichlet}%
\end{equation}
where $\partial w_{j}(v)$ denotes the derivative of $w$ at the vertex $v$
taken along the edge $e_{j}$ in the direction outward the vertex. The sum in
the Kirchhoff-Neumann condition (\ref{KN}) is taken over all the edges $e_{j}%
$, $j=1,\ldots,M$.

Thus, if for a certain value of the spectral parameter $\lambda$ there exists
a nontrivial solution $w$ of the problem (\ref{Schr})-(\ref{Dirichlet}), it
consists of $M$ components $w_{j}(x)$ associated with the corresponding edges
$e_{j}$, such that each component satisfies the equation
\begin{equation}
-w_{i}^{\prime\prime}(x)+q_{i}(x)w_{i}(x)=\lambda w_{i}(x),\quad x\in
(0,L_{i}), \label{Schri}%
\end{equation}
the initial condition
\begin{equation}
w_{i}(0)=0 \label{Dirichleti}%
\end{equation}
and the coupling relations
\begin{equation}
w_{i}(L_{i})=w_{j}(L_{j})\quad\text{for all }i,j\in\left\{  1,\ldots
,M\right\}  , \label{continuity}%
\end{equation}
\begin{equation}
\sum_{i=1}^{M}w_{i}^{\prime}(L_{i})=0. \label{KNi}%
\end{equation}
The potential $q_{i}(x)$ in (\ref{Schri}) is the component of $q(x)$
corresponding to the edge $e_{i}$.

It is well known that problem (\ref{Schr})-(\ref{Dirichlet}) has a discrete
spectrum consisting of an infinite sequence of eigenvalues $\lambda_{1}%
\leq\lambda_{2}\leq\ldots$, $\lambda_{k}\rightarrow+\infty$, and the
corresponding eigenfunctions $u_{1}$, $u_{2}$,... can be chosen so that
$\left\{  u_{k}\right\}  _{k=1}^{\infty}$ forms an orthonormal basis in
$\mathcal{H}:=L_{2}(\Omega)$,%
\[
\left\langle u_{i},u_{j}\right\rangle _{\mathcal{H}}=\int_{\Omega}%
u_{i}(x)u_{j}(x)dx=\delta_{ij}.
\]

Denote $\chi_{k}(\gamma):=\partial u_{k}(\gamma)$, $\gamma\in\Gamma$. Let
$\alpha_{k}$ be the $M$-component column vector defined as
\[
\alpha_{k}=\operatorname{col}\left(  \frac{\chi_{k}(\gamma)}{\sqrt{\lambda
_{k}}}\right)  _{\gamma\in\Gamma}.
\]

\begin{definition}
The set of pairs $\left\{  \lambda_{k},\alpha_{k}\right\}  _{k=1}^{\infty}$ is
called the Dirichlet spectral data of the graph $\Omega$.
\end{definition}

The inverse spectral problem consists in recovering the potential $q(x)$ from
the Dirichlet spectral data. The aim of the present work is to propose a
method for its solution.

\section{Neumann series of Bessel functions representations}

By $\varphi_{i}(\rho,x)$ and $S_{i}(\rho,x)$ we denote the solutions of the
equation
\[
-y^{\prime\prime}(x)+q_{i}(x)y(x)=\rho^{2}y(x),\quad x\in(0,L_{i})
\]
satisfying the initial conditions
\[
\varphi_{i}(\rho,0)=1,\quad\varphi_{i}^{\prime}(\rho,0)=0,
\]%
\[
S_{i}(\rho,0)=0,\quad S_{i}^{\prime}(\rho,0)=1.
\]
Here $\rho=\sqrt{\lambda}$, $\operatorname{Im}\rho\geq0$. The main tool used
in the present work are the series representations obtained in \cite{KNT} for
the solutions of (\ref{Schri}) and their derivatives.

\begin{theorem}
[\cite{KNT}]\label{Th NSBF} The solutions $\varphi_{i}(\rho,x)$ and
$S_{i}(\rho,x)$ and their derivatives with respect to $x$ admit the following
series representations
\begin{align}
\varphi_{i}(\rho,x)  &  =\cos\left(  \rho x\right)  +\sum_{n=0}^{\infty
}(-1)^{n}g_{i,n}(x)\mathbf{j}_{2n}(\rho x),\label{phiNSBF}\\
S_{i}(\rho,x)  &  =\frac{\sin\left(  \rho x\right)  }{\rho}+\frac{1}{\rho}%
\sum_{n=0}^{\infty}(-1)^{n}s_{i,n}(x)\mathbf{j}_{2n+1}(\rho x),\label{S}\\
\varphi_{i}^{\prime}(\rho,x)  &  =-\rho\sin\left(  \rho x\right)  +\frac{1}%
{2}\int_{0}^{x}q_{i}(t)\,dt\,\cos\left(  \rho x\right)  +\sum_{n=0}^{\infty
}(-1)^{n}\gamma_{i,n}(x)\mathbf{j}_{2n}(\rho x),\label{phiprimeNSBF}\\
S_{i}^{\prime}(\rho,x)  &  =\cos\left(  \rho x\right)  +\frac{1}{2\rho}%
\int_{0}^{x}q_{i}(t)\,dt\,\sin\left(  \rho x\right)  +\frac{1}{\rho}\sum
_{n=0}^{\infty}(-1)^{n}\sigma_{i,n}(x)\mathbf{j}_{2n+1}(\rho x),
\label{Sprime}%
\end{align}
where $\mathbf{j}_{k}(z)$ stands for the spherical Bessel function of order
$k$ (see, e.g., \cite{AbramowitzStegunSpF}). The coefficients $g_{i,n}(x)$,
$s_{i,n}(x)$, $\gamma_{i,n}(x)$ and $\sigma_{i,n}(x)$ can be calculated
following a simple recurrent integration procedure (see \cite{KNT} or
\cite[Sect. 9.4]{KrBook2020}), starting with
\begin{align}
g_{i,0}(x)  &  =\varphi_{i}(0,x)-1,\quad s_{i,0}(x)=3\left(  \frac{S_{i}%
(0,x)}{x}-1\right)  ,\label{beta0}\\
\gamma_{i,0}(x)  &  =g_{i,0}^{\prime}(x)-\frac{1}{2}\int_{0}^{x}%
q_{i}(t)\,dt,\quad\sigma_{i,0}(x)=\frac{s_{i,0}(x)}{x}+s_{i,0}^{\prime
}(x)-\frac{3}{2}\int_{0}^{x}q_{i}(t)\,dt.\nonumber
\end{align}
For every $\rho\in\mathbb{C}$ all the series converge pointwise. For every
$x\in\left[  0,L_{i}\right]  $ the series converge uniformly on any compact
set of the complex plane of the variable $\rho$, and the remainders of their
partial sums admit estimates independent of $\operatorname{Re}\rho$.
\end{theorem}

This last feature of the series representations (the independence of
$\operatorname{Re}\rho$ of the estimates for the remainders) is of crucial
importance for what follows. In particular, it means that for $S_{i,N}%
(\rho,x):=\frac{\sin\left(  \rho x\right)  }{\rho}+\frac{1}{\rho}\sum
_{n=0}^{N}(-1)^{n}s_{i,n}(x)\mathbf{j}_{2n+1}(\rho x)$ (and analogously for
$S_{i,N}^{\prime}(\rho,x)$) the estimate holds
\begin{equation}
\left\vert S_{i}(\rho,x)-S_{i,N}(\rho,x)\right\vert <\varepsilon_{i,N}(x)
\label{estim S}%
\end{equation}
for all $\rho\in\mathbb{R}$, where $\varepsilon_{i,N}(x)$ is a positive
function tending to zero as $N\rightarrow\infty$. That is, the approximate
solution $S_{i,N}(\rho,x)$ approximates the exact one equally well for small
and for large values of $\rho$. This is especially convenient when considering
direct and inverse spectral problems. Moreover, for a fixed $z$ the numbers
$\mathbf{j}_{k}(z)$ rapidly decrease as $k\rightarrow\infty$, see, e.g.,
\cite[(9.1.62)]{AbramowitzStegunSpF}. Hence, the convergence rate of the
series for any fixed $\rho$ is, in fact, exponential.

More detailed estimates for the series remainders depending on the regularity
of the potential can be found in \cite{KNT}.

Note that formulas (\ref{beta0}) indicate that the potential $q_{i}(x)$ can be
recovered from the first coefficients of the series (\ref{phiNSBF}) or
(\ref{S}). We have
\begin{equation}
q_{i}(x)=\frac{g_{i,0}^{\prime\prime}(x)}{g_{i,0}(x)+1} \label{qi from g0}%
\end{equation}
and
\begin{equation}
q_{i}(x)=\frac{\left(  xs_{i,0}(x)\right)  ^{\prime\prime}}{xs_{i,0}(x)+3x}.
\label{qi from s0}%
\end{equation}

\section{Solution of the direct problem\label{Sect Direct}}

Due to the Dirichlet condition at the boundary vertices the component $w_{i}$
of the eigenfunction has the form%
\[
w_{i}(x)=c_{i}(\rho)S_{i}(\rho,x),
\]
where $c_{i}(\rho)$ is a constant. Equalities (\ref{continuity}) give us $M-1$
equations of the form%
\[
c_{i}(\rho)S_{i}(\rho,L_{i})-c_{i+1}(\rho)S_{i+1}(\rho,L_{i+1})=0,\quad
i=1,\ldots,M-1,
\]
and equality (\ref{KNi}) gives us the $M$-th equation%
\[
\sum_{i=1}^{M}c_{i}(\rho)S_{i}^{\prime}(\rho,L_{i})=0.
\]
This system of $M$ equations can be written in the matrix form%
\[
\mathbf{S}(\rho)\overrightarrow{c}(\rho)=
\]%
\[
\left(
\begin{array}
[c]{cccccc}%
S_{1}(\rho,L_{1}) & -S_{2}(\rho,L_{2}) & 0 & 0 & \ldots & 0\\
0 & S_{2}(\rho,L_{2}) & -S_{3}(\rho,L_{3}) & 0 & \ldots & 0\\
&  & \ddots &  &  & \\
0 & 0 & \ldots & 0 & S_{M-1}(\rho,L_{M-1}) & -S_{M}(\rho,L_{M})\\
S_{1}^{\prime}(\rho,L_{1}) & S_{2}^{\prime}(\rho,L_{2}) & \ldots &  &
S_{M-1}^{\prime}(\rho,L_{M-1}) & S_{M}^{\prime}(\rho,L_{M})
\end{array}
\right)  \left(
\begin{array}
[c]{c}%
c_{1}(\rho)\\
c_{2}(\rho)\\
\vdots\\
c_{M-1}(\rho)\\
c_{M}(\rho)
\end{array}
\right)  =0.
\]
The number $\rho$ is a square root of an eigenvalue if and only if the
determinant of the matrix $\mathbf{S}(\rho)$ is zero. Computation of zeros of
$\det\mathbf{S}(\rho)$ can be easily performed by using the series
representations (\ref{S}) and (\ref{Sprime}). Let $\rho_{k}$ be such a zero.
The vector $\overrightarrow{c}(\rho_{k})$ then is an eigenvector of
$\mathbf{S}(\rho_{k})$ corresponding to the zero eigenvalue of $\mathbf{S}%
(\rho_{k})$. In order to compute the column vector $\alpha_{k}$ we notice that
$\alpha_{k}=a\overrightarrow{c}(\rho_{k})$, where $a$ is a constant such that
\[
\sum_{i=1}^{M}\left(  ac_{i}(\rho_{k})\right)  ^{2}\rho_{k}^{2}\int_{0}%
^{L_{i}}S_{i}^{2}(\rho_{k},x)dx=1,
\]
which means that the corresponding eigenfunction is normalized. Thus,%
\[
a^{2}=\left(  \rho_{k}^{2}\sum_{i=1}^{M}c_{i}^{2}(\rho_{k})\int_{0}^{L_{i}%
}S_{i}^{2}(\rho_{k},x)dx\right)  ^{-1}%
\]
and hence the components of $\alpha_{k}$ have the form
\begin{equation}
\alpha_{k,i}=\frac{c_{i}(\rho_{k})}{\rho_{k}\sqrt{\sum_{i=1}^{M}c_{i}^{2}%
(\rho_{k})\int_{0}^{L_{i}}S_{i}^{2}(\rho_{k},x)dx}}. \label{alphaki}%
\end{equation}
The components of the normalized eigenfunction $u_{k}$ have the form
$u_{k,i}(x)=\alpha_{k,i}\rho_{k}S_{i}(\rho_{k},x)$.

\section{Solution of the inverse spectral problem\label{Sect Inverse}}

The method for solving the inverse problem consists of several steps. Given
the spectral data, first, we compute the constants $\left\{  s_{i,n}%
(L_{i})\right\}  $ and $\left\{  \sigma_{i,n}(L_{i})\right\}  $ for every
$i=1,\ldots,M$, that is, the values of the coefficients from (\ref{S}) and
(\ref{Sprime}) at the endpoint of the edge $e_{i}$. This first step allows us
to split the problem on the graph and reduce it to the problems on the edges.
This is because the knowledge of the numbers $\left\{  s_{i,n}(L_{i})\right\}
$ and $\left\{  \sigma_{i,n}(L_{i})\right\}  $ allows us in the second step to
compute the Dirichlet-Dirichlet and Dirichlet-Neumann spectra for the
potential $q_{i}(x)$, $x\in\left[  0,L_{i}\right]  $, thus obtaining a
two-spectra inverse problem for each component of the potential $q(x)$.

In the third step the two-spectra inverse problem is solved with the aid of
the representation (\ref{S}) and an analogous series representation for the
solution $\psi_{i}(\rho,x)$ of (\ref{Schri}) satisfying the initial conditions%
\begin{equation}
\psi_{i}(\rho,L_{i})=1,\quad\psi_{i}^{\prime}(\rho,L_{i})=0. \label{psi init}%
\end{equation}
The two-spectra inverse problem reduces to a system of linear algebraic
equations, from which we obtain the coefficient $s_{i,0}(x)$. Finally, the
potential $q_{i}(x)$ is calculated from (\ref{qi from s0}).

\subsection{Reduction to the two spectra inverse problem on the edge
\label{Subsect reduction to two spectra}}

Given the set of the Dirichlet spectral data, that is, given the numbers
$\rho_{k}$ and the constants $c_{i}(\rho_{k})$ for which the equalities are
valid
\begin{equation}
c_{i}(\rho_{k})S_{i}(\rho_{k},L_{i})-c_{i+1}(\rho_{k})S_{i+1}(\rho_{k}%
,L_{i+1})=0,\quad i=1,\ldots,M-1 \label{cont rhok}%
\end{equation}
and
\begin{equation}
\sum_{i=1}^{M}c_{i}(\rho_{k})S_{i}^{\prime}(\rho_{k},L_{i})=0. \label{KN rhok}%
\end{equation}
Note that $c_{i}(\rho_{k})=\chi_{k}(\gamma_{i})=\alpha_{k,i}\rho_{k}$.

Consider equations (\ref{cont rhok}). With the aid of (\ref{S}) they can be
written in the form
\begin{align}
&  c_{i}(\rho_{k})\sum_{n=0}^{\infty}(-1)^{n}s_{i,n}(L_{i})\mathbf{j}%
_{2n+1}(\rho_{k}L_{i})-c_{i+1}(\rho_{k})\sum_{n=0}^{\infty}(-1)^{n}%
s_{i+1,n}(L_{i+1})\mathbf{j}_{2n+1}(\rho_{k}L_{i+1})\nonumber\\
&  =-c_{i}(\rho_{k})\sin(\rho_{k}L_{i})+c_{i+1}(\rho_{k})\sin(\rho_{k}%
L_{i+1}),\quad i=1,\ldots,M-1. \label{eqSi}%
\end{align}

The set of these equations gives us a system of linear algebraic equations for
computing the constants $s_{i,n}(L_{i})$, $i=1,\ldots,M$. For every $\rho_{k}$
we have $M-1$ equations of the form (\ref{eqSi}). The knowledge of a finite
number of spectral data allows us to compute a finite number of the
approximate constants $\left\{  s_{i,n}(L_{i})\right\}  _{n=0}^{N}$ for all
$i=1,\ldots,M$. Thus, for each edge $e_{i}$ we can compute the approximate
function
\begin{equation}
S_{i,N}(\rho,L_{i})=\frac{\sin\left(  \rho L_{i}\right)  }{\rho}+\frac{1}%
{\rho}\sum_{n=0}^{N}(-1)^{n}s_{i,n}(L_{i})\mathbf{j}_{2n+1}(\rho L_{i})
\label{SiN}%
\end{equation}
for any value of $\rho$. Note that small perturbations in the coefficients
$s_{i,n}(L_{i})$, caused by truncating series in (\ref{eqSi}) and disposing of
a finite number of spectral data, provoke only small perturbations in the
value of $S_{i,N}(\rho,L_{i})$ due to the fact that the spherical Bessel
functions are bounded ($\left\vert \mathbf{j}_{m}(t)\right\vert \leq1$) and
decreasing ($\mathbf{j}_{m}(t)=\frac{\sin(t-m\pi/2)}{t}+O\left(  \frac
{1}{t^{2}}\right)  $) for all $m=0,1,\ldots$ and $t\in\mathbb{R}$.

The function (\ref{SiN}) gives us an approximation of the solution $S_{i}%
(\rho,x)$ evaluated at the endpoint $L_{i}$. Moreover, for all real $\rho$
(and $\rho$ remaining in some strip $\left\vert \operatorname{Im}%
\rho\right\vert <C$) the difference between the exact value $S_{i}(\rho
,L_{i})$ and the approximate one $S_{i,N}(\rho,L_{i})$ remains bounded by the
same constant $\varepsilon_{i,N}(L_{i})$ (see (\ref{estim S})).

Since $S_{i}(\rho,x)$ satisfies the Dirichlet condition at $x=0$, zeros of the
function $S_{i}(\rho,L_{i})$, which we denote by $\mu_{i,k}$, are the square
roots of the eigenvalues of the Sturm-Liouville problem%
\begin{equation}
-y^{\prime\prime}+q_{i}(x)y=\lambda y,\quad x\in(0,L_{i}), \label{SLi}%
\end{equation}%
\begin{equation}
y(0)=y(L_{i})=0. \label{DD cond}%
\end{equation}
Thus, computation of zeros of the function $S_{i,N}(\rho,L_{i})$ gives us an
approximation of the Dirichlet-Dirichlet spectrum of the potential $q_{i}(x)$
(the spectrum of problem (\ref{SLi}), (\ref{DD cond})). The uniform bound
(\ref{estim S}) guarantees the uniform bound of the absolute error of
computation of the numbers $\mu_{i,k}$ \cite[Proposition 7.1]{KT2015JCAM}.
Moreover, arbitrarily many numbers $\mu_{i,k}$ can be computed (hundreds or
even thousands in practice). This gives us the possibility to compute an
important parameter
\[
\omega_{i}=\frac{1}{2}\int_{0}^{L_{i}}q_{i}(t)\,dt.
\]
Indeed, the sequence of numbers $\mu_{i,k}$ satisfies the asymptotic relation%
\[
\mu_{i,k}=\frac{\pi}{L_{i}}k+\frac{\omega_{i}}{\pi k}+\frac{\xi_{i,k}}{k},
\]
where $\left\{  \xi_{i,k}\right\}  _{k=1}^{\infty}\in\ell_{2}$ (see, e.g.,
\cite[p. 18]{Yurko2007}). In other words, the sequence $\left\{  k\left(
\mu_{i,k}-\frac{\pi}{L_{i}}k\right)  -\frac{\omega_{i}}{\pi}\right\}
_{k=1}^{\infty}$ belongs to $\ell_{2}$. Suppose that finitely many numbers
$\left\{  \mu_{i,k}\right\}  _{k=1}^{K_{D}}$ are known. Then an approximate
value of the parameter $\omega_{i}$ can be found by minimizing the $\ell_{2}%
$-norm of the sequence $\left\{  k\left(  \mu_{i,k}-\frac{\pi}{L_{i}}k\right)
-\frac{\omega_{i}}{\pi}\right\}  _{k=K_{0}}^{K_{D}}$, where $K_{0}$ is some
integer between $1$ and $K_{D}$, chosen to skip several first numbers
$\mu_{i,k}$ which can differ considerably from their asymptotics. One can
take, e.g., the floor of $K_{D}/2$, $K_{0}=\left\lfloor K_{D}/2\right\rfloor
$, and obtain that
\[
\omega_{i}\approx\arg\min_{\omega_{i}}\sum_{k=\left\lfloor K_{D}%
/2\right\rfloor }^{K_{D}}\left(  k\left(  \mu_{i,k}-\frac{\pi}{L_{i}}k\right)
-\frac{\omega_{i}}{\pi}\right)  ^{2}.
\]

Having computed the parameter $\omega_{i}$, with the aid of (\ref{Sprime})
equation (\ref{KN rhok}) can be written in the form%
\begin{align*}
&  c_{1}(\rho_{k})\sum_{n=0}^{\infty}(-1)^{n}\sigma_{1,n}(L_{1})\mathbf{j}%
_{2n+1}(\rho_{k}L_{1})+\ldots+c_{M}(\rho_{k})\sum_{n=0}^{\infty}(-1)^{n}%
\sigma_{M,n}(L_{M})\mathbf{j}_{2n+1}(\rho_{k}L_{M})\\
&  =-\sum_{i=1}^{M}c_{i}(\rho_{k})\left(  \rho_{k}\cos\left(  \rho_{k}%
L_{i}\right)  -\omega_{i}\sin\left(  \rho_{k}L_{i}\right)  \right)  .
\end{align*}
This gives us a system of linear algebraic equations for computing the
constants $\sigma_{i,n}(L_{i})$.

Thus, for each edge $e_{i}$ we can compute the value of the function
\begin{equation}
S_{i,N}^{\prime}(\rho,L_{i})=\cos\left(  \rho L_{i}\right)  +\frac{\omega_{i}%
}{\rho}\sin\left(  \rho L_{i}\right)  +\frac{1}{\rho}\sum_{n=0}^{N}%
(-1)^{n}\sigma_{i,n}(L_{i})\mathbf{j}_{2n+1}(\rho L_{i}) \label{SiNprime}%
\end{equation}
for any value of $\rho$, which gives us the approximate derivative of the
corresponding solution evaluated at the endpoint $L_{i}$, and again,
$\left\vert S_{i}^{\prime}(\rho,L_{i})-S_{i,N}^{\prime}(\rho,L_{i})\right\vert
$ remains uniformly bounded for all $\rho$ belonging to a strip $\left\vert
\operatorname{Im}\rho\right\vert <C$.

With the aid of $S_{i,N}^{\prime}(\rho,L_{i})$ we compute the approximate
Dirichlet-Neumann eigenvalues of the potential $q_{i}(x)$. Indeed, considering
the Sturm-Liouville problem
\begin{equation}
-y^{\prime\prime}+q_{i}(x)y=\lambda y,\quad x\in(0,L_{i}), \label{SLi2}%
\end{equation}%
\begin{equation}
y(0)=y^{\prime}(L_{i})=0 \label{DN cond}%
\end{equation}
we notice that zeros of the function $S_{i}^{\prime}(\rho,L_{i})$, which we
denote by $\nu_{i,k}$, are the square roots of the eigenvalues of problem
(\ref{SLi2}), (\ref{DN cond}). Thus, having computed a number of the constants
$\left\{  \sigma_{i,n}(L_{i})\right\}  _{n=0}^{N}$ gives us the possibility to
compute the approximate Dirichlet-Neumann singular numbers $\nu_{i,k}$ by
computing zeros of the function (\ref{SiNprime}). Again, arbitrarily many
numbers $\nu_{i,k}$ can be computed, and thus on each edge $e_{i}$ we obtain a
classical inverse Sturm-Liouville problem of recovering the potential
$q_{i}(x)$ from two spectra: the Dirichlet-Dirichlet spectrum of problem
(\ref{SLi}), (\ref{DD cond}) and Dirichlet-Neumann spectrum of problem
(\ref{SLi2}), (\ref{DN cond}).

\subsection{Solution of the inverse spectral problem on the edge}

At this stage we dispose of two sequences of singular numbers $\left\{
\mu_{i,k}\right\}  _{k=1}^{K_{D}}$ and $\left\{  \nu_{i,k}\right\}
_{k=1}^{K_{N}}$ which are square roots of the eigenvalues of problems
(\ref{SLi}), (\ref{DD cond}) and (\ref{SLi2}), (\ref{DN cond}), respectively,
as well as of two sequences of numbers $\left\{  s_{i,n}(L_{i})\right\}
_{n=0}^{N}$ and $\left\{  \sigma_{i,n}(L_{i})\right\}  _{n=0}^{N}$, which are
the values of the coefficients from (\ref{S}) and (\ref{Sprime}) at the
endpoint. The number $K_{N}$ may be chosen equal to $K_{D}$ or different.

The two-spectra inverse problem, in principle, can be solved by any available
method. However, we chose an approach which makes use of the knowledge of the
coefficients $\left\{  s_{i,n}(L_{i})\right\}  _{n=0}^{N}$ and $\left\{
\sigma_{i,n}(L_{i})\right\}  _{n=0}^{N}$ and does not require the
Dirichlet-Dirichlet singular numbers $\left\{  \mu_{i,k}\right\}
_{k=1}^{K_{D}}$ (which were used in the previous step for computing the
parameters $\omega_{i}$).

Let us consider the solution $\psi_{i}(\rho,x)$ of equation (\ref{SLi2})
satisfying the initial conditions at $L_{i}$:
\[
\psi_{i}(\rho,L_{i})=1,\quad\psi_{i}^{\prime}(\rho,L_{i})=0.
\]
Analogously to the solution (\ref{phiNSBF}) the solution $\psi_{i}(\rho,x)$
admits the series representation
\begin{equation}
\psi_{i}(\rho,x)=\cos\left(  \rho\left(  L_{i}-x\right)  \right)  +\sum
_{n=0}^{\infty}\left(  -1\right)  ^{n}\tau_{i,n}\left(  x\right)
\mathbf{j}_{2n}\left(  \rho\left(  L_{i}-x\right)  \right)  , \label{psi}%
\end{equation}
where $\tau_{i,n}\left(  x\right)  $ are corresponding coefficients, analogous
to $g_{i,n}\left(  x\right)  $ from (\ref{phiNSBF}). Similarly to
(\ref{qi from g0}) the equality
\begin{equation}
q_{i}(x)=\frac{\tau_{i,0}^{\prime\prime}(x)}{\tau_{i,0}(x)+1}
\label{q from tau}%
\end{equation}
is valid.

Note that for $\rho=\nu_{i,k}$ the solutions $S_{i}(\nu_{i,k},x)$ and
$\psi_{i}(\nu_{i,k},x)$ are linearly dependent because both are eigenfunctions
of problem (\ref{SLi2}), (\ref{DN cond}). Hence there exist such real
constants $\beta_{i,k}\neq0$, that
\begin{equation}
S_{i}(\nu_{i,k},x)=\beta_{i,k}\psi_{i}(\nu_{i,k},x). \label{S=psi}%
\end{equation}
Moreover, these multiplier constants can be easily calculated by recalling
that $\psi_{i}(\nu_{i,k},L_{i})=1$. Thus,%
\begin{align}
\beta_{i,k}  &  =S_{i}(\nu_{i,k},L_{i})\approx S_{i,N}(\nu_{i,k}%
,L_{i})\nonumber\\
&  =\frac{\sin\left(  \nu_{i,k}L_{i}\right)  }{\nu_{i,k}}+\frac{1}{\nu_{i,k}%
}\sum_{n=0}^{N}(-1)^{n}s_{i,n}(L_{i})\mathbf{j}_{2n+1}(\nu_{i,k}L_{i}).
\label{betaik}%
\end{align}
Having computed these constants we use equation (\ref{S=psi}) for constructing
a system of linear algebraic equations for the coefficients $s_{i,n}(x)$ and
$\tau_{i,n}\left(  x\right)  $. Indeed, equation (\ref{S=psi}) can be written
in the form%
\begin{align*}
&  \frac{1}{\nu_{i,k}}\sum_{n=0}^{\infty}(-1)^{n}s_{i,n}(x)\mathbf{j}%
_{2n+1}(\nu_{i,k}x)-\beta_{i,k}\sum_{n=0}^{\infty}\left(  -1\right)  ^{n}%
\tau_{i,n}\left(  x\right)  \mathbf{j}_{2n}\left(  \nu_{i,k}\left(
L_{i}-x\right)  \right) \\
&  =-\frac{\sin\left(  \nu_{i,k}x\right)  }{\nu_{i,k}}+\beta_{i,k}\cos\left(
\nu_{i,k}\left(  L_{i}-x\right)  \right)  .
\end{align*}

We have as many of such equations as many Dirichlet-Neumann singular numbers
$\nu_{i,k}$ are computed. For computational purposes we choose some natural
number $N_{c}$ - the number of the coefficients $s_{i,n}(x)$ and $\tau
_{i,n}\left(  x\right)  $ to be computed. More precisely, we choose a
sufficiently dense set of points $x_{m}\in(0,L_{i})$ and at every $x_{m}$
consider the equations
\begin{align}
&  \frac{1}{\nu_{i,k}}\sum_{n=0}^{N_{c}}(-1)^{n}s_{i,n}(x_{m})\mathbf{j}%
_{2n+1}(\nu_{i,k}x_{m})-\beta_{i,k}\sum_{n=0}^{N_{c}}\left(  -1\right)
^{n}\tau_{i,n}\left(  x_{m}\right)  \mathbf{j}_{2n}\left(  \nu_{i,k}\left(
L_{i}-x_{m}\right)  \right) \nonumber\\
&  =-\frac{\sin\left(  \nu_{i,k}x_{m}\right)  }{\nu_{i,k}}+\beta_{i,k}%
\cos\left(  \nu_{i,k}\left(  L_{i}-x_{m}\right)  \right)  ,\quad
k=1,\ldots,K_{N}. \label{last eq}%
\end{align}
Solving this system of equations we find $s_{i,0}(x_{m})$ and $\tau
_{i,0}\left(  x_{m}\right)  $ and consequently $s_{i,0}(x)$ and $\tau
_{i,0}\left(  x\right)  $ at a dense set of points of the interval $(0,L_{i}%
)$. Finally, with the aid of (\ref{qi from s0}) or (\ref{q from tau}) we
compute $q_{i}(x)$.

Schematically the proposed method for solving the inverse spectral problem on
a quantum star graph is presented in the following diagram.%
\[
SD=\left\{  \rho_{k},\alpha_{k}\right\}  _{k=1}^{K}\quad\overset
{(1)}{\Longrightarrow}\quad\left\{  s_{i,n}(L_{i}),\omega_{i};\ \sigma
_{i,n}(L_{i})\right\}  _{n=0}^{N}\quad\overset{(2)}{\Longrightarrow}%
\quad\text{two spectra }\left\{  \mu_{i,k}\right\}  _{k=1}^{K_{D}},\ \left\{
\nu_{i,k}\right\}  _{k=1}^{K_{N}}%
\]%
\[
\left\{  s_{i,n}(L_{i})\right\}  _{n=0}^{N},\,\left\{  \nu_{i,k}\right\}
_{k=1}^{K_{N}}\quad\overset{(3)}{\Longrightarrow}\quad\left\{  \beta
_{i,k}\right\}  _{k=1}^{K_{N}}%
\]%
\[
\left\{  \nu_{i,k},\beta_{i,k}\right\}  _{k=1}^{K_{N}}\quad\overset
{(4)}{\Longrightarrow}\quad s_{i,0}(x),\ \tau_{i,0}\left(  x\right)
\quad\overset{(5)}{\Longrightarrow}\quad q_{i}(x).
\]

Note that in the first step $K$ should be greater than $N$.

After step (1) the problem is reduced to separate problems on the edges. The
two spectra problem arising after step (2) can be solved by different existing
methods, nevertheless here we propose a method which uses the fact that the
constants $\left\{  s_{i,n}(L_{i})\right\}  _{n=0}^{N}$ and $\left\{
\sigma_{i,n}(L_{i})\right\}  _{n=0}^{N}$ are also known.

It is worth mentioning that when dealing with the truncated systems of linear
algebraic equations we do not seek to work with the square systems. In
computations a least-squares solution of an overdetermined system (provided by
Matlab, which we used in this work) gives very satisfactory results.

\section{Numerical examples}

\textbf{Example 1. }Let us consider a star graph of five edges of lengths
\begin{equation}
L_{1}=\frac{e}{2},\,L_{2}=1,\,L_{3}=\frac{\pi}{2},\,L_{4}=\frac{\pi}%
{3},\,L_{5}=\frac{e^{2}}{4}.\, \label{L}%
\end{equation}
The corresponding five components of the potential are defined as follows%
\begin{equation}
q_{1}(x)=\left\vert x-1\right\vert +1,\,q_{2}(x)=e^{-(x-\frac{1}{2})^{2}%
},\,q_{3}(x)=\sin\left(  8x\right)  +\frac{2\pi}{3}, \label{q1}%
\end{equation}%
\begin{equation}
q_{4}(x)=\cos\left(  9x^{2}\right)  +1,\quad q_{5}(x)=\frac{1}{x+0.1}.
\label{q2}%
\end{equation}

We computed $100$ eigenvalues\ and corresponding norming constants as
described in Section \ref{Sect Direct}. In Table 1 some of the computed
eigenvalues are presented.

\bigskip%
\begin{tabular}
[c]{|l|l|l|l|l|}\hline
\multicolumn{5}{|l|}{Table 1: Computed eigenvalues of the quantum star graph
(\ref{L})-(\ref{q2})}\\\hline
$\lambda_{1}$ & $\lambda_{2}$ & $\lambda_{5}$ & $\lambda_{10}$ &
$\lambda_{100}$\\\hline
$1.5656490615325$ & $2.1509437903100$ & $3.2180647998489$ & $5.493521290269$ &
$46.683594634217$\\\hline
\end{tabular}

\medskip

Then to the computed spectral data we applied the method introduced in Section
\ref{Sect Inverse} for solving the inverse problem. The number $N$ of the
coefficients in (\ref{SiN}) and (\ref{SiNprime}) was chosen\ $N=10$. Our
numerical experiments showed that even a more reduced number of the
coefficients (say, five-six of them) gives already quite satisfactory
(slightly less accurate) results. For the number $N_{c}$ from (\ref{last eq})
was chosen the same value, $N_{c}=10$, though this parameter, of course, may
differ from $N$.

In Fig. 1 we show the recovered potentials on the five edges. The absolute
error of each recovered potential slightly grows at the endpoints of the
interval attaining its maximum value\ in the second decimal digit (for
$q_{5}(x)$ at the origin), while in the interior points of the interval it
remains two orders of magnitude smaller.%

\begin{figure}
[ptb]
\begin{center}
\includegraphics[
height=2.5482in,
width=5.5534in
]%
{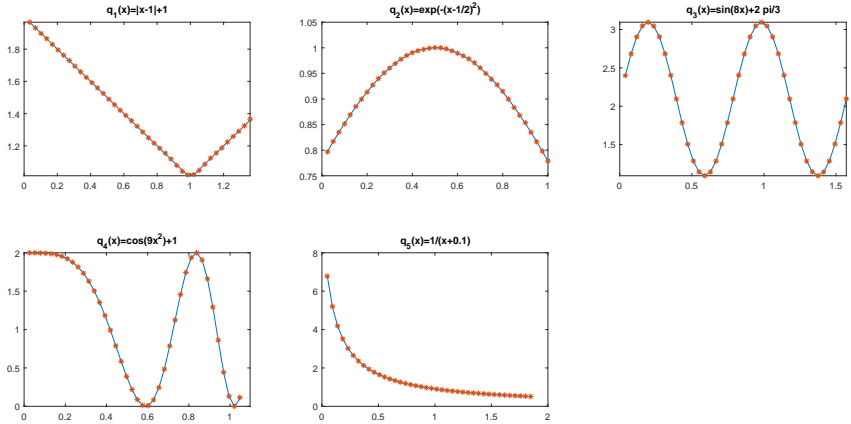}%
\caption{The recovered potential of the quantum star graph (\ref{L}%
)-(\ref{q2}) from $100$ spectral eigendata with $N=10$.}%
\label{Fig1}%
\end{center}
\end{figure}

It is interesting to track how accurately the two spectra were computed on
each edge. For example, Table 2 presents some of the \textquotedblleft
exact\textquotedblright\ Dirichlet-Dirichlet eigenvalues on $e_{2}$ computed
with the aid of the Matslise package (first column), the approximate
eigenvalues, computed as described in Subsection
\ref{Subsect reduction to two spectra} by calculating zeros of $S_{2,10}%
(\rho,L_{2})$ and the absolute error of each presented eigenvalue. Notice that
both the absolute and the relative errors remain small even for large indices.

\medskip

\bigskip%
\begin{tabular}
[c]{|l|l|l|l|}\hline
\multicolumn{4}{|l|}{Table 2: Dirichlet-Dirichlet eigenvalues of $q_{2}(x)$%
}\\\hline
$n$ & $\lambda_{n}$ & $\widetilde{\lambda}_{n}$ & $\left\vert \lambda
_{n}-\widetilde{\lambda}_{n}\right\vert $\\\hline
$1$ & $10.8381543$ & $10.8381555$ & $1.2\cdot10^{-6}$\\\hline
$11$ & $1195.1450218$ & $1195.1450282$ & $6.4\cdot10^{-6}$\\\hline
$41$ & $16591.72758$ & $16591.72741$ & $0.00016$\\\hline
$101$ & $100680.75706$ & $100680.75692$ & $0.00014$\\\hline
$201$ & $398742.80997$ & $398742.80983$ & $0.00013$\\\hline
\end{tabular}

\bigskip

Based on the computed Dirichlet-Dirichlet eigenvalues the parameter
$\omega_{i}$ was computed for all $i$, as explained in Subsection
\ref{Subsect reduction to two spectra}. In the case of $e_{2}$ the absolute
error of the approximation was $\left\vert \omega_{2}-\widetilde{\omega}%
_{2}\right\vert =7.2\cdot10^{-5}$.

Next, using the quantum graph spectral eigendata and the computed $\omega_{i}$
the Dirichlet-Neumann eigenvalues for each potential component $q_{i}(x)$ were
computed. Table 3 presents the result of computation for $q_{2}(x)$. One can
observe that while the accuracy is slightly worse than in the case of the
Dirichlet-Dirichlet eigenvalues, it is still quite satisfactory and uniform.

\bigskip%

\begin{tabular}
[c]{|l|l|l|l|}\hline
\multicolumn{4}{|l|}{Table 3: Dirichlet-Neumann eigenvalues of $q_{2}(x)$%
}\\\hline
$n$ & $\lambda_{n}$ & $\widetilde{\lambda}_{n}$ & $\left\vert \lambda
_{n}-\widetilde{\lambda}_{n}\right\vert $\\\hline
$1$ & $3.3898$ & $3.3920$ & $0.0022$\\\hline
$11$ & $1089.0464$ & $1089.0487$ & $0.0023$\\\hline
$41$ & $16189.5411$ & $16189.5418$ & $0.0006$\\\hline
$101$ & $99686.394414$ & $99686.394415$ & $8.8\cdot10^{-7}$\\\hline
$201$ & $396761.48688$ & $396761.48678$ & $0.0001$\\\hline
\end{tabular}

\medskip

\bigskip

In the considered example the worst accuracy was obtained for the component of
the potential $q_{5}(x)$. Table 4 shows the accuracy of the
Dirichlet-Dirichlet eigenvalues computed.\medskip

\bigskip%

\begin{tabular}
[c]{|l|l|l|l|}\hline
\multicolumn{4}{|l|}{Table 4: Dirichlet-Dirichlet eigenvalues of $q_{5}(x)$%
}\\\hline
$n$ & $\lambda_{n}$ & $\widetilde{\lambda}_{n}$ & $\left\vert \lambda
_{n}-\widetilde{\lambda}_{n}\right\vert $\\\hline
$1$ & $3.99363768$ & $3.99363773$ & $4.5\cdot10^{-8}$\\\hline
$11$ & $351.545497$ & $351.545501$ & $3.8\cdot10^{-6}$\\\hline
$41$ & $4863.543766$ & $4863.543738$ & $2.8\cdot10^{-5}$\\\hline
$101$ & $29505.854685672381$ & $29505.8546235738$ & $6.2\cdot10^{-5}$\\\hline
\end{tabular}

\medskip

\bigskip

Based on this result, the approximation of the parameter $\omega_{5}$ resulted
in the following value of the absolute error $\left\vert \omega_{5}%
-\widetilde{\omega}_{5}\right\vert =0.0042$. Finally, some Dirichlet-Neumann
eigenvalues for $q_{5}(x)$ are presented in Table 5.

\medskip%

\begin{tabular}
[c]{|l|l|l|l|}\hline
\multicolumn{4}{|l|}{Table 5: Dirichlet-Neumann eigenvalues of $q_{5}(x)$%
}\\\hline
$n$ & $\lambda_{n}$ & $\widetilde{\lambda}_{n}$ & $\left\vert \lambda
_{n}-\widetilde{\lambda}_{n}\right\vert $\\\hline
$1$ & $1.5067$ & $1.5026$ & $0.0040$\\\hline
$11$ & $320.4508$ & $320.4470$ & $0.0037$\\\hline
$41$ & $4745.682$ & $4745.678$ & $0.0041$\\\hline
$101$ & $29214.456$ & $29214.452$ & $0.0044$\\\hline
\end{tabular}

\medskip

\bigskip

It is worth noting that the whole computation takes few seconds performed in
Matlab 2017 on a Laptop equipped with a Core i7 Intel processor.

The method copes equally well with inverse spectral problems on star graphs
with a larger number of edges, though to obtain similar accuracy in this case
more spectral data are required.

\textbf{Example 2.} Consider a nine edges quantum star graph with five edges
and components of potential coinciding with those from the previous example
and with the additional edges $L_{6}=1.1,\,L_{7}=1.2,\,L_{8}=1.3,\,L_{9}=1.4$
and the corresponding components of the potential
\[
q_{6}(x)=\frac{1}{\left(  x+0.1\right)  ^{2}},\,\,q_{7}(x)=e^{x}%
,\,\,q_{8}(x)=\pi^{2},\,\,q_{9}(x)=J_{0}(9x),
\]
where $J_{0}(z)$ stands for the Bessel function of the first kind of order
zero. Fig. 2 presents the potential recovered from $200$ spectral data.%

\begin{figure}
[ptb]
\begin{center}
\includegraphics[
height=3.3715in,
width=4.4895in
]%
{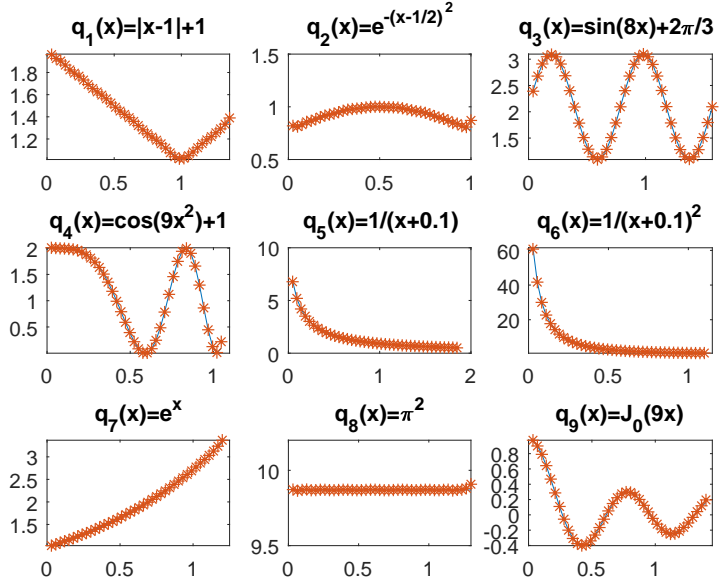}%
\caption{The potential of the quantum star graph from Example 2 recovered from
$200$ spectral eigendata with $N=10$.}%
\label{Fig2}%
\end{center}
\end{figure}

\section{Conclusions}

A new method for solving the inverse spectral problem on quantum star graphs
is developed. The main role in the proposed approach is played by the
coefficients of the Neumann series of Bessel functions expansion of solutions
of the Sturm-Liouville equation. With their aid the given spectral data lead
to separate two spectra Sturm-Liouville inverse problems on each edge. These
two-spectra problems are solved by a direct method reducing each problem to a
system of linear algebraic equations, and the crucial observation is that the
potential is recovered from the first component of the solution vector.

The method is simple, direct and accurate. Its performance is illustrated by
numerical examples. In subsequent works we plan to extend this method to other
types of inverse spectral problems and to more general graphs.

\section*{Acknowledgements}

Research was supported by CONACYT, Mexico via the project 284470 and partially
performed at the Regional mathematical center of the Southern Federal
University with the support of the Ministry of Science and Higher Education of
Russia, agreement 075-02-2022-893. The research of Sergei Avdonin was
supported in part by the National Science Foundation, grant DMS 1909869, and
by Moscow Center for Fundamental and Applied Mathematics.


\begin{thebibliography}{99}                                                                                               %


\bibitem {AbramowitzStegunSpF}Abramovitz M. and Stegun I. A. (1972),
\textit{Handbook of mathematical functions}, New York: Dover.

\bibitem {ArioliBenzi2018}Arioli M. and Benzi M. (2018), A finite element
method for quantum graphs. \textit{IMA J. Numer. Anal.}, \textbf{38}, no. 3, 1119--1163.

\bibitem {AvdoninBelinskiyMatthews}{Avdonin S., Belinskiy B., and Matthews J.}
(2011), {Inverse problem on the semi-axis: local approach}, \textit{Tamkang
Journal of Mathematics}, \textbf{42}, no. 3, 1--19.

\bibitem {AvdoninBell2015}Avdonin S. and Bell J. (2015), Determining physical
parameters for a neuronal cable model defined on a tree graph, \textit{Journal
of Inverse Problems and Imaging,} \textbf{9}, no. 3, 645-659.

\bibitem {AvdoninKravchenko2022}Avdonin S. and Kravchenko V. V. (2022) Method
for solving inverse spectral problems on quantum star graphs. Submitted.

\bibitem {AvdoninKurasov2008}Avdonin S. and Kurasov P. (2008), {Inverse
problems for quantum trees,} \textit{Inverse Problems and Imaging},
\textbf{2}, no. 1, 1--21.

\bibitem {AvdoninLeugeringMikhaylov2010}Avdonin S., Leugering G., and
Mikhaylov V. (2010), {On an inverse problem for tree-like networks of elastic
strings}, \textit{Zeit. Angew. Math. Mech.}, \textbf{90}, no. 2, 136--150.

\bibitem {Belishev Vakulenko 2006}Belishev M. and Vakulenko A. (2006), Inverse
problems on graphs: Recovering the tree of strings by the BC-method, \emph{J.
Inv. Ill-Posed Problems}, \textbf{14} , 29-46.

\bibitem {BerkolaikoKuchment}Berkolaiko G. and Kuchment P. (2013),
\textit{Introduction to Quantum Graphs}, AMS, Providence, R.I.

\bibitem {Brown et al 2003}Brown B. M., Samko V. S., Knowles I. W., Marletta
M. (2003), Inverse spectral problem for the Sturm--Liouville equation,
\emph{Inverse Probl.} \textbf{19} , 235--252.

\bibitem {Chadan et al 1997}Chadan Kh., Colton D., P\"{a}iv\"{a}rinta L.,
Rundell W. (1997), \emph{An introduction to inverse scattering and inverse
spectral problems}. SIAM, Philadelphia.

\bibitem {Kr2019JIIP}Kravchenko V. V. (2019), On a method for solving the
inverse Sturm--Liouville problem,\textbf{ }\emph{J. Inverse Ill-posed Probl.
}\textbf{27}, 401--407.

\bibitem {KrBook2020}Kravchenko V. V. (2020), \emph{Direct and inverse
Sturm-Liouville problems: A method of solution}, Birkh\"{a}user, Cham.

\bibitem {KNT}Kravchenko V. V., Navarro L. J. and Torba\ S. M. (2017),
Representation of solutions to the one-dimensional Schr\"{o}dinger equation in
terms of Neumann series of Bessel functions, \emph{Appl. Math. Comput.}
\textbf{314}, 173--192.

\bibitem {KT2015JCAM}Kravchenko V. V. and Torba\ S. M. (2015), Analytic
approximation of transmutation operators and applications to highly accurate
solution of spectral problems,\emph{ Journal of Computational and Applied
Mathematics }\textbf{275}, 1-26.

\bibitem {KT2021 IP1}Kravchenko V. V. and Torba S. M. (2021), A direct method
for solving inverse Sturm-Liouville problems, \emph{Inverse Probl.
}\textbf{37}, 015015 (32pp).

\bibitem {Kurasov}Kurasov P. (2022) Quantum Graphs: Spectral Theory and
Inverse Problems, Springer (to appear).

\bibitem {Kurasov et al 2005}Kurasov P. and Nowaczyk M. (2005), Inverse
spectral problem for quantum graphs, \emph{J. Phys. A.,} \textbf{38}, 4901-4915.

\bibitem {LagneseLeugeringSchmidt1994}Lagnese J.E., Leugering G., and Schmidt
E.J.P.G. (1994), \emph{ Modelling, Analysis and Control of Multi-Link Flexible
Structures}, Basel: Birkh\"{a}user.

\bibitem {LevitanInverse}Levitan B. M. (1987), \emph{Inverse Sturm-Liouville
problems}, VSP, Zeist.

\bibitem {Mugnolo}Mugnolo, D. (2014), \textit{Semigroup Methods for Evolution
Equations on Networks,} Understanding Complex Systems, Springer, Cham.

\bibitem {SitnikShishkina Elsevier}Shishkina E. L. and Sitnik S. M. (2020),
\emph{Transmutations, singular and fractional differential equations with
applications to mathematical physics}, Elsevier, Amsterdam.

\bibitem {SavchukShkalikov}Savchuk A. M., Shkalikov A. A. (2005), Inverse
problem for Sturm--Liouville operators with distribution potentials:
reconstruction from two spectra, \emph{Russ. J. Math. Phys. }\textbf{12}, 507--514.

\bibitem {Yurko2005}Yurko V. A. (2005), Inverse Sturm-Lioville operator on
graphs, \emph{Inverse Problems}, \textbf{21}, 1075-1086.

\bibitem {Yurko2007}Yurko V. A. (2007), \emph{Introduction to the theory of
inverse spectral problems}, Fizmatlit, Moscow, (in Russian).
\end{thebibliography}
\end{document}